\documentclass{amsart}
\usepackage{amsmath,amsthm,amssymb}
\usepackage{hyperref}
\usepackage{tikz,color}

\usepackage{a4wide}

\title{Proof of a conjecture of M\'esz\'aros and Morales on the
  volume of a flow polytope}

\author{Jang Soo Kim}
\address[Jang Soo Kim]{
Department of Mathematics, Sungkyunkwan University, Suwon 440-746,
South Korea}
\email{jangsookim@skku.edu}
\date{\today}

\newtheorem{thm}{Theorem}

\theoremstyle{definition}

\newtheorem{conj}{Conjecture}

\newcommand\Cat{\operatorname{Cat}}
\newcommand\CT{\operatorname{CT}}

\begin{document}

\maketitle

\begin{abstract}
  We prove a conjecture of M\'esz\'aros and Morales on the volume of a
  flow polytope. Independently from our work, Zeilberger sketched a
  proof of their conjecture. In fact, our proof is the same as
  Zeilberger's proof. The purpose of this note is to give a more
  detailed proof of the conjecture.
\end{abstract}

\section{Introduction}

We prove a conjecture of M\'esz\'aros and Morales \cite{MM} on the
volume of a flow polytope, which is a type $D$ analog of the
Chan-Robbins-Yuen polytope \cite{CRY2000}. Independently from our
work, Zeilberger also proved the conjecture and sketched his proof in
\cite{Z}. In fact, our proof is the same as Zeilberger's proof. The
purpose of this note is to give a more detailed proof of the
conjecture.

For a function $f(z)$ with a Laurent series expansion at $z$, we
denote by $\CT_z f(z)$ the constant term of the Laurent expansion of
$f(z)$ at $0$. In other words, if $f(z)=\sum_{n=-\infty}^\infty a_n
z^n$, then $\CT_z f(z) = a_0$.

The conjecture of M\'esz\'aros and Morales can be stated as the
following constant term identity.

\begin{conj} \cite[Conjecture~7.6]{MM} \label{MMconj}
For an integer $n\ge2$, we have
\[
  \CT_{x_n}\CT_{x_{n-1}}\cdots \CT_{x_1}
\prod_{j=1}^n x_j^{-1}(1-x_j)^{-2}
\prod_{1\le j<k\le  n} (x_k-x_j)^{-1} (1-x_k-x_j)^{-1} 
= 2^{n^2} \prod_{k=1}^n \Cat(k),
\]  
where $\Cat(k)=\frac1{k+1}\binom{2k}{k}$. 
\end{conj}

Conjecture~\ref{MMconj} is a type $D_n$ analog of the constant term
identity
\begin{equation}
  \label{eq:1}
  \CT_{x_n}\CT_{x_{n-1}}\cdots \CT_{x_1}
\prod_{j=1}^n (1-x_j)^{-2}
\prod_{1\le j<k\le  n} (x_k-x_j)^{-1} 
=\prod_{k=1}^n \Cat(k),
\end{equation}
which was conjectured by Chan, Robbins, and Yuen \cite{CRY2000} and
proved by Zeilberger \cite{Zeilberger1999}.

In Section~\ref{sec:zeilb-proof-chan} we recall Zeilberger's
reformulation of Morris' constant term identity. In
Section~\ref{sec:proof-conjecture} we prove Conjecture~\ref{MMconj}.

\section{Zeilberger's reformulation of Morris' constant term identity}
\label{sec:zeilb-proof-chan}

Zeilberger \cite{Zeilberger1999} rewrote Morris' identity
\cite{Morris} as follows: For nonnegative integers $a$ and $b$ and a
positive half integer $c$, we have
\begin{equation}
  \label{eq:zeilberger}
\CT_{x_n}\CT_{x_{n-1}}\cdots \CT_{x_1} \prod_{i=1}^n (1-x_i)^{-a} x_i^{-b} 
\prod_{1\le i<j\le  n}(x_j-x_i)^{-2c} 
= \frac{1}{n!} \prod_{j=0}^{n-1} \frac{\Gamma(a+b+(n-1+j)c)\Gamma(c)}
{\Gamma(a+jc)\Gamma(c+jc)\Gamma(b+jc+1)}.
\end{equation}
Zeilberger \cite{Zeilberger1999} showed that, when $a=2,b=0,c=1/2$,
the right hand side of \eqref{eq:zeilberger} is
\begin{equation}
  \label{eq:cat}
\frac{1}{n!} \prod_{j=0}^{n-1}
\frac{\Gamma(\frac{n+3+j}2)\Gamma(\frac 12)}
{\Gamma(\frac{4+j}2)\Gamma(\frac{1+j}2)\Gamma(\frac{2+j}2)}=  
\prod_{k=1}^n \Cat(k),
\end{equation}
which together with \eqref{eq:zeilberger} implies \eqref{eq:1}. 

By Cauchy's integral formula, if $f(z)$ has a Laurent series expansion
at $0$, we have
\begin{equation}
  \label{eq:CT}
\CT_z f(z) = \frac{1}{2\pi i} \oint_C \frac{f(z)} z dz,
\end{equation}
where $C$ is the circle $\{z: |z|=\epsilon\}$ oriented
counterclockwise for a real number $\epsilon >0$ such that $f(z)$ is
holomorphic inside $C$ except $0$.  Thus \eqref{eq:zeilberger} can be rewritten
as
\begin{multline}
  \label{eq:zeilberger2}
\frac{1}{(2\pi i)^n}\oint_{C_n}\dots\oint_{C_1}  \prod_{j=1}^n (1-x_j)^{-a} x_j^{-b-1} 
\prod_{1\le j<k\le  n}(x_k-x_j)^{-2c}  dx_1\cdots dx_n\\
= \frac{1}{n!} \prod_{j=0}^{n-1} \frac{\Gamma(a+b+(n-1+j)c)\Gamma(c)}
{\Gamma(a+jc)\Gamma(c+jc)\Gamma(b+jc+1)},
\end{multline}
where $C_j$ is the circle $\{z: |z|=j\epsilon\}$ oriented
counterclockwise for a number $0<\epsilon<\frac 1n$.

\section{Proof of the conjecture of M\'esz\'aros and Morales}
\label{sec:proof-conjecture}

We will prove the following theorem. 
\begin{thm}\label{thm} For a nonnegative integer $a$ and a positive half integer
  $c$, we have
\begin{multline*}
\CT_{x_n}\CT_{x_{n-1}}\cdots \CT_{x_1}  \prod_{j=1}^n x_j^{-a+1}(1-x_j)^{-a}
\prod_{1\le j<k\le  n} (x_j-x_k)^{-2c} (1-x_j-x_k)^{-2c}  
\\= 2^{2an+4c\binom n2-2n} 
\frac{1}{n!} \prod_{j=0}^{n-1} \frac{\Gamma(a-\frac12+(n-1+j)c)\Gamma(c)}
{\Gamma(\frac12+jc)\Gamma(c+jc)\Gamma(a+jc)}.
\end{multline*}
\end{thm}

We can get Conjecture~\ref{MMconj} as the special case $a=2,c=1/2$ of
Theorem~\ref{thm} as follows. When $a=2,c=1/2$, the right hand side of
the formula in Theorem~\ref{thm} is
\[
\frac{2^{n^2+n}}{n!}\prod_{j=0}^{n-1}
\frac{\Gamma(\frac{n+2+j}2)\Gamma(\frac 12)}
{\Gamma(\frac{1+j}2)\Gamma(\frac{1+j}2)\Gamma(\frac{4+j}2)}.
\]
Since
\[
\prod_{j=0}^{n-1}
\frac{\Gamma(\frac{n+3+j}2)\Gamma(\frac 12)}
{\Gamma(\frac{4+j}2)\Gamma(\frac{1+j}2)\Gamma(\frac{2+j}2)}
\left/
\prod_{j=0}^{n-1}
\frac{\Gamma(\frac{n+2+j}2)\Gamma(\frac 12)}
{\Gamma(\frac{1+j}2)\Gamma(\frac{1+j}2)\Gamma(\frac{4+j}2)}=  
\right.
\frac{\Gamma(\frac{2n+2}2)\Gamma(\frac{1}2)}
{\Gamma(\frac{n+2}2)\Gamma(\frac{n+1}2)}
=2^n, 
\]
we obtain Conjecture~\ref{MMconj} by \eqref{eq:cat}.

For the rest of this section we prove Theorem~\ref{thm}.  From now on
we assume that $\epsilon>0$ is a very small number.

By \eqref{eq:CT} we have
\begin{align*}
&  \CT_{x_n}\CT_{x_{n-1}}\cdots \CT_{x_1}  \prod_{j=1}^n x_j^{-a+1}(1-x_j)^{-a}
\prod_{1\le j<k\le  n} (x_k-x_j)^{-2c} (1-x_k-x_j)^{-2c} \\
&= \frac{1}{(2\pi i)^n}\oint_{C_n}\cdots\oint_{C_1}
\prod_{j=1}^n x_j^{-a}(1-x_j)^{-a}
\prod_{1\le j<k\le  n} (x_k-x_j)^{-2c} (1-x_k-x_j)^{-2c}
dx_1\cdots dx_n,
\end{align*}
where $C_j$ is the circle $\{z: |z|=j\epsilon\}$ oriented
counterclockwise. 

Using the change of variables $x_j = \frac{1-z_j}2$ or $z_j=1-2x_j$,
the above is equal to
\begin{multline*}
\frac{1}{(2\pi i)^n}\oint_{C'_n}\cdots\oint_{C'_1}
\prod_{j=1}^n \left(\frac{1-z_j}2\right)^{-a}
\left(\frac{1+z_j}2\right)^{-a} \\
\times \prod_{1\le j<k\le  n} \left(\frac{-z_k+z_j}2\right)^{-2c}
\left(\frac{z_k+z_j}2\right)^{-2c} (-2)^{-n} dz_1\cdots dz_n\\
=\frac{(-1)^n2^{2an+4c\binom n2-n}}{(2\pi i)^n}\oint_{C'_n}\cdots\oint_{C'_1}
\prod_{j=1}^n (1-z_j^2)^{-a}
\prod_{1\le j<k\le  n} (z_j^2-z_k^2)^{-2c}
dz_1\cdots dz_n,
\end{multline*}
where $C'_j$ is the circle $\{z: |z-1|=2j\epsilon\}$ oriented
counterclockwise.

Using the change of variables $z_j^2 = y_j$ or $z_j=y_j^{1/2}$, the above is
equal to
\begin{align*}
&\frac{(-1)^n2^{2an+4c\binom n2-2n}}{(2\pi i)^n}\oint_{C''_n}\cdots\oint_{C''_1}
\prod_{j=1}^n (1-y_j)^{-a} y_j^{-1/2}
\prod_{1\le j<k\le  n} (y_j-y_k)^{-2c}
dy_1\cdots dy_n,
\end{align*}
where $C''_j$ is the circle $\{z: |z-1|=4j\epsilon\}$ oriented
counterclockwise. This is because if $C_j'$ is parametrized by
$1+2j\epsilon e^{i\theta}$ for $0\le \theta\le 2\pi$, then the image
of $C_j'$ under the map $z\mapsto z^2$ can be parametrized by
$1+4j\epsilon e^{i\theta} + 4j^2\epsilon^2 e^{2i\theta}$ for $0\le
\theta\le 2\pi$. Since $\epsilon$ is very small we can deform this
image to the circle $C_j''$ without changing the contour integral.

Using the change of variables $t_j = 1-y_j$, the above is
equal to
\begin{align*}
\frac{2^{2an+4c\binom n2-2n}}{(2\pi i)^n} \oint_{C'''_n}\cdots\oint_{C'''_1}
\prod_{j=1}^n t_j^{-a} (1-t_j)^{-1/2}
\prod_{1\le j<k\le  n} (t_k-t_j)^{-2c}
dt_1\cdots dt_n,
\end{align*}
where $C'''_j$ is the circle $\{z: |z|=4j\epsilon\}$ oriented
counterclockwise.  

Using \eqref{eq:zeilberger2}, we finish the proof of
Theorem~\ref{thm}. 

\section*{Acknowledgement}
The author would like to thank Dennis Stanton for helpful discussion.

\end{document}